# PHILIP HALL'S PROBLEM ON NON-ABELIAN SPLITTERS

Rüdiger Göbel and Saharon Shelah

**Abstract**

Philip Hall raised around 1965 the following question which is stated in the Kourovka Notebook [11, p. 88]: *Is there a non-trivial group which is isomorphic with every proper extension of itself by itself?* We will decompose the problem into two parts: We want to find non-commutative splitters, that are groups $G \neq 1$ with $\mathrm{Ext}\,(G,G) = 1$. The class of splitters fortunately is quite large so that extra properties can be added to $G$. We can consider groups $G$ with the following properties: There is a complete group $L$ with cartesian product $L^\omega \cong G$, $\mathrm{Hom}\,(L^\omega, S_\omega) = 0$ ($S_\omega$ the infinite symmetric group acting on $\omega$) and $\mathrm{End}\,(L,L) = \mathrm{Inn}\,L \cup \{0\}$. We will show that these properties ensure that $G$ is a splitter and hence obviously a Hall-group in the above sense. Then we will apply a recent result from our joint paper [8] which also shows that such groups exist, in fact there is a class of Hall-groups which is not a set.

## 1 Introduction

Philip Hall investigated in his lectures at Cambridge around 1965 the following class of groups, which are characterized by our first

**Definition 1.1** *We will say that a group $G$ is a Hall-group if any extension $H$ of $G$ is isomorphic to $G$ provided $G$ normal in $H$ and $H/G \cong G$.*

John Lennox communicated in Kourovka's Notebook Hall's question concerning the existence of these Hall-groups. We want to show the existence of Hall-groups using

1991 *Mathematics Subject Classification:* Primary 13C05, 18E40, 18G05, 20K20, 20K35, 20K40; Secondary: 13D30, 18G25, 20K25, 20K30, 13C10
The authors are supported by the project No. G 0545-173, 06/97 of the German-Israeli Foundation for Scientific Research & Development.
[GbSh:738] in Shelah's list of publications.

some terminology which turned out to be particular important in module theory (and abelian groups) recently, see [9, 10].

**Definition 1.2** *A group $G$ is a weak splitter if any extension of $G$ by $G$ splits. In homological terms this is to say that $\mathrm{Ext}\,(G,G) = 1$, or equivalently any short exact sequence*

$$1 \longrightarrow G \xrightarrow{\beta} H \xrightarrow{\alpha} G \to 1$$

*gives raise to a splitting map $\gamma : G \longrightarrow H$ such that $\gamma\alpha = id_G$. Here maps are acting on the right. Hence $H = \ker \alpha \rtimes \mathrm{Im}\,\gamma \cong G \rtimes \mathrm{Im}\,\gamma$ and if $\mathrm{Im}\,\gamma$ is also normal in $H$, then we say that $G$ is a splitter. Hence $G \triangleleft H$ with $H/G \cong G$ implies $H = G \times N$ for some normal subgroup $N \triangleleft H$ which is isomorphic to $G$.*

Recall that $G = N \rtimes U$ is the semidirect product of $N$ and $U$, where $N$ is normal in $G$, and if also $U$ is normal, then we write $N \times U$ for the direct product.

In case of abelian groups classical splitters are free abelian groups as well as torsion-free cotorsion groups, which are well-known for a long time, see Fuchs [4]. Other splitters were constructed only recently, see [9, 10]. They were also fundamental for solving the flat cover conjecture for modules. Here we will study non-commutative splitters, which however are obtained quite differently. Such groups will be based on the following *ad hoc*

**Definition 1.3** *We will say that a group $L \neq 1$ is rigid if the following holds:*

(i) *$L$ is complete, i.e. $L$ has trivial center $\mathfrak{z}L = 1$ and no outer automorphisms $\mathrm{Aut}\, L = \mathrm{Inn}\, L$.*

(ii) *$\mathrm{Hom}\,(L^\omega, S_\omega) = 0$.*

(iii) *$\mathrm{End}\, L = \mathrm{Inn}\, L \cup \{0\}$.*

Here $S_\omega$ is the full symmetric group acting on a countable set $\omega = \{0, 1, 2, \ldots\}$. $\mathrm{Aut}\, L$ and $\mathrm{Inn}\, L$ denote the automorphism group and the group of inner automorphisms of $L$, respectively. Moreover $\mathrm{Hom}\,(A, B)$ is the set of homomorphisms from the group $A$ to the group $B$, where $0$ is the zero-homomorphism mapping any element to $1$; in particular $\mathrm{End}\, A = \mathrm{Hom}\,(A, A)$ is the near endomorphism ring of $A$. By $A^\kappa$ we denote the cartesian power over the cardinal $\kappa$ of the group $A$. Using these natural definitions we have a possibility to find Hall-group by means of our

**Main Theorem 1.4** *If $L$ is a rigid group and $G \cong L^\omega$, then the following holds*

(a) $\mathrm{Aut}\, G = \mathrm{Inn}\, G \rtimes S_\omega$.

(b) $G$ is a splitter.

(c) $G \times G \cong G$.

Here we note that many groups with $(a)$ in the Main Theorem 1.4 are constructed in [6, 5, 3, 2] - in fact for arbitrary groups in place of $S_\omega$, but $(b), (c)$ are also crucial. We have an immediate

**Corollary 1.5** *If $L$ is a rigid group, then $L^\omega$ is a Hall-group.*

In view of Theorem 1.4 it remains show the existence of a rigid group. In Section 3 where we will provide the following result from [8] and sketch the essentials of its (lengthy) proof in a few lines. We also note that condition $(i)$ of Definition 1.3 does not speak about the group $L$ but only of its infinite cartesian power. It would be desirable for convenience and for esthetic reasons to have a group theoretic condition immediately on $L$. This is established in Proposition 2.1 and leads to our main theorem from [8] mentioned as Theorem 3.1 in the last section. Here are the parts needed for application in the proof of our Main Theorem 1.4.

**Theorem 1.6** *For any infinite cardinals $\kappa < \mu$ with $\kappa$ regular, $\mu = \mu^\kappa$ and $\lambda = \mu^+ > 2^{\aleph_0}$ there is a group $H$ with the following properties.*

  (i) *$H$ is a simple group of cardinality $\lambda$.*

 (ii) *There are an element $h \in H$ such that any element of $H$ is a product of at most 4 conjugates of $h$.*

(iii) *$H$ is rigid.*

## 2 Proof of Main Theorem 1.4

PROOF. Condition $(c)$ is obvious because $G \times G = L^\omega \times L^\omega \cong L^\omega = G$.
In order to show $(a)$ let $S_\omega$ act naturally on $\omega$ and write $G = \prod_{n \in \omega} Le_n$, hence

$$x = \prod_{n \in \omega} x_n e_n \text{ with } x_n \in L \tag{2.1}$$

denotes a general element of $G$. Also let $[x] = \{n \in \omega : x_n \neq 1\}$ denote the support of $x$. Moreover

$$G_A = \{x \in G : [x] \subseteq A\} \subseteq G \text{ for any } A \subseteq \omega.$$

If $\pi \in S_\omega$, then $\pi$ induces an automorphism of $G$ (also denoted by $\pi$) given by
$$\pi : G \longrightarrow G \ (x = \prod_{n \in \omega} x_n e_n \longrightarrow x\pi = \prod_{n \in \omega} x_{n\pi} e_n).$$
Hence $S_\omega \subseteq \operatorname{Aut} G$ and if $g \in G$, then let
$$g^* : G \longrightarrow G \ (x \longrightarrow xg^* = g^{-1}xg)$$
denote the inner automorphism, conjugation by $g$. Hence
$$G^* = \operatorname{Inn} G = \{g^* : g \in G\}$$
is normal in $\operatorname{Aut} G$ and visibly $G^* \cap S_\omega = 1$. The semidirect product
$$G^* \rtimes S_\omega \subseteq \operatorname{Aut} G$$
is a subgroup of $\operatorname{Aut} G$ and we claim that the two groups are equal. Let $\sigma \in \operatorname{Aut} G$ be a given automorphism and $n, m \in \omega$. Then we define a homomorphism
$$\sigma_{nm} : L \longrightarrow L \ (x \longrightarrow xe_n \longrightarrow xe_n\sigma \longrightarrow (xe_n\sigma)_m),$$
where clearly $(xe_n\sigma)_m$ is the $m$–th coordinate of $xe_n\sigma = \prod_{i \in \omega}(xe_n\sigma)_i e_i$ as follows from our convention (2.1). Using the canonical embedding
$$\iota_n : L \longrightarrow Le_n \subseteq G \ (x \longrightarrow xe_n)$$
and the canonical projection
$$\pi_m : G \longrightarrow L \ (x = \prod_{n \in \omega} x_n e_n \longrightarrow x_m),$$
we can also say that
$$\sigma_{nm} = \iota_n \sigma \pi_m \in \operatorname{End} L.$$
From Definition 1.3 (iii) we have $\sigma_{mn} \in \operatorname{Inn} L \cup \{0\}$ and if $\sigma_{mn} \neq 0$ then there is $t_{\sigma mn} \in L$ such that $\sigma_{mn} = t^*_{\sigma mn}$. Let
$$Y_\sigma = \{(n, m) : \exists \ t^*_{\sigma mn} = \sigma_{mn}\},$$
hence
$$\sigma_{mn} = 0 \iff (n, m) \in \omega \times \omega \setminus Y_\sigma.$$
Now we define $A_n = \{m \in \omega : (n, m) \in Y_\sigma\}$ and claim that

$$\text{the } A_n \ (n \in \omega) \text{ are pair-wise disjoint.} \tag{2.2}$$

If $m \in A_{n_1} \cap A_{n_2}$ and $n_1 \neq n_2$ then $(n_1, m), (n_2, m) \in Y_\sigma$. Take any $x_i \in L$ and let $y_i = x_i e_{n_i}$ ($i = 1, 2$). Hence $[y_1, y_2] = 1$ in $G$ and $[L\sigma_{n_1 m}, L \sigma_{n_2 m}] = 1$ follows from $L\sigma_{n_i m} = L e_{n_i} \sigma \pi_m \subseteq L$. If $L\sigma_{n_i m} = L$, then $[L, L] = 1$ and $L$ would be abelian, which contradicts $\mathfrak{z} L = 1 \neq L$. If $L\sigma_{n_1 m} \neq L$ then $\sigma_{n_1 m} = 0$ by Definition 1.3 $(iii)$, which contradicts $(n_1, m) \in Y_\sigma$ and (2.2) is shown.

Next we observe:

$$\text{If } x = \prod_{i \in \omega} x_i e_i \in G \text{ and } x_n = 1 \text{ then } x\sigma\pi_m = 1 \text{ for all } m \in A_n. \tag{2.3}$$

Note that $x \in G_{\omega \setminus \{n\}}$, $G = G_{\{n\}} \times G_{\omega \setminus \{n\}}$ and $x\sigma\pi_m \in G_{\omega \setminus \{n\}} \sigma \pi_m$. From

$$[G_{\{n\}}, G_{\omega \setminus \{n\}}] = 1$$

follows that $G_{\omega \setminus \{n\}} \sigma \pi_m$ and $G_{\{n\}} \sigma \pi_m$ commute, and $m \in A_n$ implies that

$$G_n \sigma \pi_m = L \sigma_{nm} = L.$$

Hence $G_{\omega \setminus \{n\}} \sigma \pi_m = 1$ from $\mathfrak{z} L = 1$ and in particular $x\sigma\pi_m = 1$.

Now we want to show that

$$|A_n| = 1 \text{ for all } n \in \omega. \tag{2.4}$$

If $A_n = \emptyset$ and $1 \neq x \in L$, then $x e_n \sigma \pi_m = 1$ for all $m \in \omega$, hence $x e_n \sigma = 1$ and $1 \neq x e_n \in \ker \sigma$ contradicts $\sigma \in \operatorname{Aut} G$. If $|A_n| > 1$, then choose $n_1 \neq n_2 \in A_n$ and define $C = G_{\omega \setminus \{n_1, n_2\}}$ and $D = G_{\{n_1, n_2\}}$. Hence $G = C \times D$, $D \cong L \times L$ and consider the canonical projection $\Pi_D : G \to D$.

From (2.3) follows $\sigma \pi_D \upharpoonright C = 0$, and $\iota_n \sigma \pi_D$ maps $L$ into $D$. This map is an isomorphism as follows from $n_1, n_2 \in A_n$, hence $L e_n \cong L e_{n_1} \oplus L e_{n_2}$. The right hand side has the outer automorphism switching coordinates, while $\operatorname{Aut} L = L^*$ by Definition 1.3$(i)$. This is impossible, and $|A_n| = 1$ follows.

In the next step we show that

$$\bigcup_{n \in \omega} A_n = \omega \tag{2.5}$$

Otherwise there is $1 \neq y \in G$ with $[y] \cap \bigcup_{n \in \omega} A_n = \emptyset$. Let $x \in G$ be the preimage of $y$ under the automorphism $\sigma$, hence $x\sigma = y$. We want to show that $x = \sum_{n \in \omega} x_n e_n = 1$, which is a contradiction. By (2.4) we can write $A_n = \{m\}$, hence

$$1 = y_m = x\sigma\pi_m = x_n \sigma_{nm} = x_n t^*_{nm}$$

using $m \notin [y]$ and (2.3) for all $n \in \omega$. We get $x_n = 1$ and $x = 1$ follows. From (2.4) and (2.5) follows that

$$\pi = \pi_\sigma : \omega \to \omega \ (n \to m) \quad \text{if } A_n = \{m\}$$

is a permutation $\pi \in S_\omega$, and it is easy to see that

$$\pi_{\sigma^{-1}} = (\pi_\sigma)^{-1} \tag{2.6}$$

We view $\pi_\sigma \in \operatorname{Aut} G$ as described at the beginning. If we replace $\sigma$ by $\sigma' = \sigma \pi_\sigma^{-1}$ then $\pi_{\sigma'} = id$ is obvious and $\sigma'$ acts component-wise, inducing conjugations $t_n^*$ for each $n \in \omega$. If $t = \prod_{n \in \omega} t_n e_n \in G$ then $\sigma' = t^*$, hence $\sigma = t^* \pi_\sigma \in G^* \rtimes S_\omega$ and (a) is shown.

It remains to show (b). Consider

$$1 \longrightarrow G \xrightarrow{id} H \xrightarrow{\beta} M \to 1 \text{ with } M \cong G$$

and let $\gamma : M \to H$ be a map of representatives for $\beta$ in $M$, that is $x\gamma\beta = x$ for all $x \in M$. If $x \in M$, then $(x\gamma)^*$ is an inner automorphism of $H$ which induces an automorphism $\alpha = (x\gamma)^* \upharpoonright G$ of $G \triangleleft H$. From (a) we have $\alpha = y_x^* \pi_x$ for some $y_x \in G$ and $\pi_x \in S_\omega$.

If we replace $\gamma$ by $\gamma' : M \longrightarrow H$ with $x\gamma' = y_x^{-1}(x\gamma)$, then $\gamma'$ is again a coset representation for $\beta$ and if we call the new map $\gamma$ again, then $y_x = 1$ for all $x \in M$.

Recall that

$$(x\gamma)^* \upharpoonright G = \pi_x \text{ for all } x \in M,$$

and consider the map

$$\pi : G \to S_\omega \ (x \to \pi_x).$$

It is easy to check that $\pi$ is a homomorphism, hence $\pi \in \operatorname{Hom}(G, S_\omega) = 0$ by Definition 1.3 (ii), which sends every $x$ to the identity in $S_\omega$. This is to say that $x\gamma \in \mathfrak{c}_H G$ for all $x \in H$. Let $C = \mathfrak{c}_H G \subseteq H$ denote the centralizer of $G$ in $H$. From $\mathfrak{z}G = 1$ follows that $G \cap C = 1$ moreover by the above $H$ is generated by the normal subgroups $G$ and $C$, hence $H = G \times C$. By the exact sequence above $\beta \upharpoonright C : C \to M \cong G$ is an isomorphism and we arrive at $H \cong G \times G$ hence $G$ is a splitter. $\square$

**Proposition 2.1** *Let $\kappa$ be a cardinal and $K$ and $L$ are groups with the following properties*

(a) $|K| < |L|$.

(b) $L$ is simple, moreover

(c) $\exists\ g \in L, m \in \omega$ such that any $x \in L$ is product of at most $m$ conjugates of $g$

Then $\mathrm{Hom}\,(L^\kappa, K) = 0$.

We derive a

**Corollary 2.2** *If $L$ is a simple group $> 2^{\aleph_0}$ such that $\mathrm{End}\,L = \mathrm{Inn}\,L \cup \{0\}$ and condition (c) of Proposition 2.1 holds, then $L$ is rigid and $L^\omega$ is a Hall-group.*

**Proof of the corollary**: We can apply Proposition 2.1 for $K = S_\omega$, hence
$$\mathrm{Hom}\,(L^\omega, S_\omega) = 0.$$
As $L$ is simple, it must be complete, so $L$ is rigid. By Corollary 1.5 it follows that $L^\omega$ is a Hall-group. $\square$

**Proof of Proposition 2.1**: Let $g \in L$ and $m \in \omega$ be as in (c). If $G = L^\kappa$ and $\sigma$ is any homomorphism from $\mathrm{Hom}\,(G, L)$, then consider any $x = \prod_{i<\kappa} e_i x_i \in G$. We want to show that $x\sigma = 1$. By (c) these coordinates $x_i \in L$ of $x$ can be expressed as
$$x_i = \prod_{j<k_i} g^{y_{ij}} \quad (i < \kappa)$$
as products of $k_i \leq m$ conjugates of $g$. If $A_k = \{i < \kappa : k_i = k\}$, then $\kappa = \bigcup_{k \leq m} A_k$ is a decomposition of $\kappa$ and we may assume $A_k = \emptyset$ for each $k \in \omega \setminus \{k_i;\ i \in \kappa\}$. Using the earlier notation we have that $G = G_{A_1} \times \cdots \times G_{A_m}$ is expressed as a direct product. If $\bar{g} = \prod_{i<\kappa} e_i g$ is the canonical diagonal element, then using the elements $y_{ij}, x, \bar{g}$ we get new elements
$$\bar{y}_{kj} = \prod_{t \in A_k} e_t y_{kj}, \quad \bar{g}_k = \prod_{t \in A_k} e_t g, \quad \bar{x}_k = \prod_{t \in A_k} e_t x_t \in G_{A_k}$$
which are restrictions of the old ones to $A_k$. In particular
$$z_{kj} = \bar{g}_k^{\bar{g}_{kj}} \in G_{A_k}, \quad \bar{x}_k = \prod_{j<k} z_{kj} \text{ and } \quad x = \prod_{k \leq m} \bar{x}_k \qquad (2.7)$$
Consider the canonical diagonal homomorphism
$$\iota_k : L \longrightarrow G_{A_k} \quad (x \longrightarrow \prod_{t \in A_k} e_t x) \qquad (k \leq m)$$
and note that $g\iota_k = \bar{g}_k$, hence
$$\bar{g}_k \sigma = g\iota_k \sigma = 1 \quad \text{from } \iota_k \sigma \in \mathrm{Hom}\,(L, K)$$
which is trivial because $|K| < |L|$ and $L$ is simple. From (2.7) follows $x\sigma = 1$ and $\mathrm{Hom}\,(G, L) = 0$. $\square$

# 3   A result from [8]

It remains to find simple groups $L$ of cardinality $|L| > 2^{\aleph_0}$ such that $\operatorname{End} L = \operatorname{Inn} L \cup \{0\}$ with the extra property that there are $g \in L, m \in \omega$ and any element of $L$ is product of at most $m$ conjugates of $g$.

Relatives of these groups are constructed in [7] in references given there and for instance in [2, 3, 5, 6]. However they are not good enough for our purpose in this paper. When working on [8] we noticed the connection to the Hall problem and extended the construction in order to incorporate its use above.

**Theorem 3.1** *Let $\mathcal{A}$ be a family of suitable groups and $\kappa < \mu$ be infinite cardinals such that $\kappa$ is regular uncountable and $\mu^\kappa = \mu$. Then we can find a group $H$ of cardinality $\lambda = \mu^+$ such that the following holds.*

   (i) *$H$ is simple. Moreover, if $1 \neq g \in H$, then any element of $H$ is a product of at most four conjugates of $g$.*

  (ii) *Any $A \in \mathcal{A}$ is a subgroup of $H$ and if $\mathcal{A}$ is not empty, then $H[\mathcal{A}] = H$, were $H[\mathcal{A}]$ is the subgroup of $H$ generated by all subgroups of $H$ isomorphic to $A$.*

 (iii) *Any monomorphism $\varphi : A \to H$ for some $A \in \mathcal{A}$ is induced by some $h \in H$, that is there is some $h \in H$ such that $\varphi = h^* \restriction A$.*

 (iv) *If $A' \subseteq H$ is an isomorphic copy of some $A \in \mathcal{A}$, then the centralizer $\mathfrak{c}_H A' = 1$ is trivial.*

  (v) *Any monomorphism $H \to H$ is an inner automorphism.*

For our application we can assume $\mathcal{A} = \emptyset$. Otherwise the following definition is needed.

**Definition 3.2** *Let $A$ be any group with trivial center and view $A \subseteq \operatorname{Aut}(A)$ as inner automorphisms of $A$. Then $A$ is called* suitable *if the following conditions hold:*

   (i) *$A \neq 1$ is a finite group.*

  (ii) *If $A' \subseteq \operatorname{Aut}(A)$ and $A' \cong A$ then $A' = A$.*

 (iii) *$\operatorname{Aut}(A)$ is complete.*

Note that $\mathrm{Aut}(A)$ has trivial center because $A$ has trivial center. Hence the last condition only requires that $\mathrm{Aut}(A)$ has no outer automorphisms. It also follows from this that any automorphism of $A$ extends to an inner automorphism of $\mathrm{Aut}(A)$. A group $A$ is complete if $A$ has trivial center $\mathfrak{z}A$ and any automorphism is inner. We also recall the easy observation from [7] which is a consequence of the classification of finite simple groups:

*All finite simple groups are suitable.*

Also note that there are many well-known examples of suitable groups which are not simple. Just apply Wieland's theorem on automorphism towers of finite groups with trivial center, see [13].

The proof of the Theorem 3.1 is a transfinite induction building the group $H$ (which has cardinal $\lambda = \mu^+$ the successor cardinal of $\mu$) as a union of a chain of subgroups $H_\alpha$ of cardinality $\mu$. The inductive steps are separated by four disjoint stationary subset $S_i$ ($0 \leq i \leq 3$) of $\lambda$, where ordinals in $S_0 \cup S_1 \cup S_2$ are limit ordinals of cofinality $\omega$ while ordinals in $S_3$ have cofinality $\kappa$. Passing from $H_\alpha$ to $H_{\alpha+1}$ now depends on the position of $\alpha$. If $\alpha$ does not belong to one of these stationary subsets, then $H_{\alpha+1} = H_\alpha * \alpha\mathbb{Z}$ is just a free product of $H_\alpha$ with a (new) infinite cyclic group $\alpha\mathbb{Z}$. If $\alpha$ belongs to one of the first three stationary subsets, then HNN-extensions are used as in [12]. In case $\alpha \in S_0$ we must deal with the conjugacy problem for $(i)$ of the theorem. Here it is enough to ensure that all elements of infinite order are conjugate and this is what HNN is designed for. Similarly we can deal with $(ii)$ and $(iii)$ by free products with amalgamated subgroup using $\alpha \in S_1$ and $\alpha \in S_2$ respectively. An enumeration of elements with repetitions ensures that nobody is forgotten. Condition $(iv)$, which is not needed her, is pure group theory, a book keeping proof by transfinite induction. The more complicated demand is Condition $(v)$ of the theorem which is a strengthening for complete. As there are obviously more possible monomorphisms to deal with then elements in the group, a combinatorial principle is needed, a Black Box must be applied which still allows to deal with the possible injections one after another while $\alpha$ runs with repetitions through $S_3$ (a set of 'only' $\lambda$ elements) while enumerating partial injective maps on the set $H$. The basic tool is that the group $H$ is build in such a way that there are many elements $h \in H_\alpha$ with large abelian centralizer $\mathfrak{c}_{H_\alpha}(h)$ of cardinality $\kappa$. These centralizer can be arranged to come from a rigid family of abelian groups, this is to say from a theorem shown more then a decade ago for abelian groups in Corner, Göbel [1, p.465]:

**Theorem 3.3** *For each subset $X \subseteq \kappa$ of the set (the cardinal) $\kappa$ there is an $\aleph_1$-free abelian group $G_X$ of cardinal $\kappa$ such that the following holds.*

$$\mathrm{Hom}(G_X, G_Y) = \begin{cases} \mathbb{Z} & : \quad \textit{if } X \subseteq Y \\ 0 & : \quad \textit{if } X \nsubseteq Y \end{cases}$$

The proof of the theorem on abelian groups uses an earlier Black Box from Shelah, see also [1] for more details.

An abelian group is $\aleph_1$-free if all its countable subgroups are free abelian. Hence many centralizers are algebraically very different. And as centralizers must be mapped under monomorphisms into centralizer, an idea often used for characterizing certain (automorphism) groups, it seems perhaps convincing that such a covering with a rigid system of abelian centralizer almost forces monomorphism to be well-behaving, to be conjugation. The details however are a bit lengthy checking done in [8].

# References

[1] A. L. S. Corner, R. Göbel, *Prescribing endomorphism algebras - a unified treatment,* Proc. London Math. Soc. (3) **50**, 447 – 479 (1985).

[2] M. Dugas, R. Göbel, *On locally finite p–groups and a problem of Philip Hall,* Journal of Algebra **159** (1993), 115 – 138

[3] M. Dugas, R. Göbel, *Torsion–free nilpotent groups and E–modules,* Archiv der Math. **54** (1990), 340 – 351

[4] L. Fuchs, *Infinite abelian groups* - Volume 1,2 Academic Press, New York 1970, 1973.

[5] R. Göbel, A. Paras, *Outer automorphism groups of countable metabelian groups,* 1998, pp. 309–317 Proceedings of the Dublin Conference on Abelian Groups, Birkhäuser Verlag, Basel 1999

[6] R. Göbel, A. Paras, *Outer automorphism groups of metabelian groups,* Journal of Pure and Appl. Algebra (2000)

[7] R. Göbel, J. L. Rodríguez and S. Shelah, *Large localizations of finite simple groups,* submitted to Crelle Journal 1999

[8] R. Göbel, S. Shelah, *Localizations of groups,* in submitted for publication

[9] R. Göbel, S. Shelah, *Cotorsion theories and splitters,* to appear in Trans. Amer. Math. Soc. (2000)

[10] R. Göbel, S. Shelah, *Almost free splitters,* Colloquium Math. **81** (1999) 193 –221.

[11] The Kourovka Notebook, *Unsolved problems in group theory,* Novosibirsk 1999

[12] R. C. Lyndon and P. E. Schupp, *Combinatorial Group Theory*, Springer Ergebnisberichte **89**, Berlin–Heidelberg–New York 1977.

[13] D. Robinson, *A Course in the Theory of Groups*, Graduate Texts in Math. vol **80**, Berlin–Heidelberg–New York 1996.

[14] P. Schultz, *Self-splitting groups,* Preprint series of the University of Western Australia at Perth (1980)

[15] P. E. Schupp, *A characterization of inner automorphisms*, Proc. Amer. Math. Soc. **101** (1987), 226–228.

Rüdiger Göbel
Fachbereich 6, Mathematik und Informatik
Universität Essen, 45117 Essen, Germany
e–mail: R.Goebel@Uni-Essen.De
and
Saharon Shelah
Department of Mathematics
Hebrew University, Jerusalem, Israel
and Rutgers University, Newbrunswick, NJ, U.S.A
e-mail: Shelah@math.huji.ac.il